\documentclass[11 pt]{amsart}
\usepackage{epsfig,graphics,amsmath,amssymb,amsthm}
\newtheorem{Theorem}{Theorem}[section]

\newtheorem{Corollary}[Theorem]{Corollary}
\newtheorem{Proposition}[Theorem]{Proposition}
\newtheorem{Lemma}[Theorem]{Lemma}

\theoremstyle{definition}

\newtheorem{Remark}[Theorem]{Remark}

\begin{document}

\title[An Algebraic Characterization of a Dehn Twist]
{An algebraic characterization of a Dehn twist for Nonorientable
Surfaces}

\author{Ferihe Atalan}
\address{Department of Mathematics, Atilim University,
06836 \newline Ankara, TURKEY} \email{ferihe.atalan@atilim.edu.tr}
\date{\today}
\thanks{The author is supported by TUBITAK-110T665.}
 \pagenumbering{arabic}

\begin{abstract} Let $N_g^k$ be a nonorientable surface of
genus \ $g\geq 5$ \ with \ $k$-punctures. In this note, we will give
an algebraic characterization of a Dehn twist about a simple closed
curve on $N_g^k$. Along the way, we will fill some little gaps in
the proofs of some theorems in \cite{A} and \cite{I1} giving
algebraic characterizations of Dehn twists about separating simple
closed curves. Indeed, our results will give an algebraic
characterization for the topological type of Dehn twists about
separating simple closed curves.
\end{abstract}

\maketitle
\section{Introduction}

In this note, $N_{g,r}^k$  will denote the  nonorientable surface of
genus $g$ with $r$ boundary components and $k$ punctures (or
distinguished points). The mapping class group of $N_{g,r}^k$, the
group of isotopy classes of all diffeomorphisms of $N_{g,r}^k$,
where diffeomorphisms and isotopies fix each point on the boundary
is denoted by ${\rm Mod}(N_{g,r}^k)$. If we restrict ourselves to
the diffeomorphisms and isotopies to those which do not permute the
punctures then we obtain the pure mapping class group ${\rm
PMod}(N_{g,r}^k)$. The subgroup ${\rm PMod^{+}}(N_{g,r}^k)$ of the
pure mapping class group consists of the pure mapping classes that
preserve the local orientation around each puncture. Also the twist
subgroup of ${\rm Mod}(N_{g,r}^k)$, generated by Dehn twists about
two-sided simple closed curves is denoted by $\mathcal{T}$.

An algebraic characterization of a Dehn twist plays important role
in the computation of the outer automorphism group of the mapping
class group of a surface (orientable or nonorientable, see \cite{A},
\cite{ASzep} and \cite{I1}). Moreover, it is one of the main tools
in the proof of the fact that any injective endomorphism of the
mapping class group of an orientable surface must be an isomorphism
proved by Ivanov and McCarthy (\cite{I3}). We note that it is also
used in the proof of the fact that any isomorphism between two
finite index subgroups of the extended mapping class group of an
orientable surface is the restriction of an inner automorphism of
this group (\cite{I2}, \cite{K1}). Using an algebraic
characterization of Dehn twists, we show that any automorphism of
the mapping class group of a surface takes Dehn twists to  Dehn
twists.

For orientable surfaces, N. V. Ivanov gave an algebraic
characterization of Dehn twists in \cite{I1}. In December 2012, E.
Irmak reported that the proofs of Ivanov's theorem on algebraic
characterization of Dehn twists about separating simple closed
curves have some gaps (see the counter example in Section\,2).

For closed nonorientable surfaces, we gave an algebraic
characterization of Dehn twists in \cite{A}, closely following
Ivanov's work \cite{I1}. Therefore, the algebraic characterization
about separating curves in \cite{A} has also gaps.  In \cite{ASzep},
using different techniques, we gave an algebraic characterization
for Dehn twists about nonseparating simple closed curves with
nonorientable complements. However, an algebraic characterization of
the Dehn twists about other two-sided simple closed curves is still
missing. In this paper, we will try to extend this result to Dehn
twists about arbitrary simple closed curves. Indeed, we will not
only give an algebraic characterization of a Dehn twist about a
simple closed curve on \ $N_g^k$ \ but also an algebraic
characterization of the topological type of the curve the Dehn twist
is about. In particular, this paper aims to fill the gaps, mentioned
above, in both \cite{A} and \cite{I1}.

The organization of the paper is as follows: In Section 2, firstly,
we will give some definitions and remanding, secondly, we will state
and prove an algebraic characterization for a power of a Dehn twist
about a simple closed curve. In Section 3, after proving some
preliminary results we will first characterize the Dehn twists about
characteristic curves on a nonorientable surface of even genus. Then
Lemma~\ref{LemChr-2.2} will lead to an algebraic characterization of
Dehn twists about separating curves (Theorem~\ref{Chr-2.1}).
Moreover, this algebraic characterization will encode the
topological type of the separating simple closed curve the Dehn
twist which is about. In Section 4, we will consider nonorientable
surfaces, where as in Section 5 we will state analogous results for
orientable surfaces.  The final section contains some immediate
consequences of the these results.

\section{Preliminaries}\label{Prelim}

Let $S$ denote the surface  $N_{g}^k$ and let $a$ be a simple closed
curve on $S$. If a regular neighborhood of $a$ is an annulus or a
M\"obius strip, then we call $a$ a two-sided or a one-sided simple
closed curve, respectively. The curve $a$ will be called trivial, if
it bounds a disc with at most one puncture or a M\"obius band on $S$
(or if it is isotopic to a boundary component). Otherwise, it is
called nontrivial.

We will denote by  $S^a$ the result of cutting of $S$ along the
simple closed curve $a$. The simple closed curve $a$ is called
nonseparating if $S^a$ is connected. Otherwise, it is called
separating.

Let $H$ be a group. If $G\leq H$ is a subgroup and $h\in H$ is an
element of $H$, then the center of $H$, the centralizer of $G$ in
$H$ and the centralizer of $h$ in $G$ will be denoted by $C(H)$,
$C_H(G)$ and $C_G(h)$, respectively.

Let \ $P: \Sigma_{g-1}^{2k}\rightarrow S$ \ be the orientation
double covering of  \ $S$ \ and \
$\tau:\Sigma_{g-1}^{2k}\rightarrow\Sigma_{g-1}^{2k}$ \ the Deck
transformation, which is an involution.  It is well known that any
diffeomorphism \ $f:S \rightarrow S$ \ has exactly two lifts to the
orientation double covering and exactly one of them is orientation
preserving. Moreover, we can regard \ ${\rm Mod}(S)$ \ as the
subgroup \ ${\rm Mod}(\Sigma_{g-1}^{2k})^\tau$, the subgroup of
mapping classes which are invariant under the action of the deck
transformation (see also \cite{A}).

Let \ $\Gamma(m)$, where \ $m \in \mathbb{Z}$, $m>1$, be the kernel
of the natural homomorphism $${\rm Mod}(\Sigma_{g-1}^{2k})
\rightarrow {\rm Aut} \,(H_{1}(\Sigma_{g-1}^{2k}, \mathbb{Z} / m
\mathbb{Z})).$$ Then \ $\Gamma(m)$ \ is a subgroup of finite index
in \ ${\rm Mod}(\Sigma_{g-1}^{2k})$. Let \ $\Gamma'(m)=\Gamma(m)
\cap {\rm Mod}(S)$, regarding \ ${\rm Mod}(S)$ \ as a subgroup of \
${\rm Mod}(\Sigma_{g-1}^{2k})$ \ as described above.

It is well known that since \ $\Gamma(m)$ consists of pure elements
only provided that \ $m\geq 3$, (\cite{I1}) and so does $\Gamma'(m)$.
In this paper, we will consider the group \ $\Gamma'(m)$ only for
$m\geq 3$. Moreover, \ $\Gamma'$ \ will always denote a subgroup of finite
index in \ $\Gamma'(m)$.

Suppose that $f \in \Gamma'(m)$  preserves a generic family of
disjoint simple closed curves $\mathcal{C}$. Then we can assume that
$f$ is identity on all points of $\mathcal{C}$.  In this case, $f$
does not interchange components of $S^{\mathcal{C}}$ and it induces
a diffeomorphism on every component of $S^{\mathcal{C}}$, which is
either isotopic to a pseudo-Anosov one, or to the identity. In this
case, $\mathcal{C}$ is called a reduction system for $f$ (and, for
its mapping class). Also, if reduction system for $f$ is minimal,
then $\mathcal{C}$ is called a minimal reduction system for $f$.

The following lemma was proved in \cite{A} mainly following Ivanov's
work.

\begin{Lemma}\label{Lem-trivial} The center of every subgroup of finite
index of $\Gamma'(m)$ is trivial.
\end{Lemma}

Let ${\rm M}(S)$ be any of three groups ${\rm PMod}^+(S)$, ${\rm
PMod}(S)$ and ${\rm Mod}(S)$.

\begin{Lemma}\label{Lem-1.1} Let $S=N_g^k$ be a connected nonorientable surface of genus $g\geq
5$ with $k$ punctures. Let $f$ be an element in $\Gamma'$. Then $f$
is a power of a Dehn twist about a simple closed curve if and only
if the following conditions are satisfied:
\begin{enumerate}
\item $C(C_{\Gamma'}(f)) \cong \mathbb{Z}$,
\item $C(C_{\Gamma'}(f)) \neq C_{\Gamma'}(f)$.
\end{enumerate}
\end{Lemma}

\begin{proof}
Assuming that the above conditions are satisfied, we need to show
that $f$ is a power of a Dehn twist about a simple closed curve.

By the second condition, $f$ cannot be pseudo-Anosov. Indeed,
if $f$ is pseudo-Anosov then the centralizer in $\Gamma'$ of a
pseudo-Anosov element is isomorphic to $\mathbb{Z}$. (We note that
the centralizer in $\Gamma$ of a pseudo-Anosov
element is isomorphic to $\mathbb{Z}$ in \cite{I1}.) Therefore,
$C(C_{\Gamma'}(f)) = C_{\Gamma'}(f)$, and we would get a contradiction.
Now using the assumptions that $C(C_{\Gamma'}(f)) \cong \mathbb{Z}$
and $f$ is a pure element, it is easy to see that $f$ is a power of a Dehn twist
about a simple closed curve.

For the other direction of the lemma, assume that $f$ is a power of
a Dehn twist about a simple closed curve $c$. To show the first
condition we will use Ivanov's ideas. We see that
$C_{\Gamma'}(f^{n})$ is equal to the set of elements of $\Gamma'$
with minimal reduction system containing $c$, up to isotopy. So,
there is a natural map from $C_{\Gamma'}(f^{n})$ to ${\rm
Mod}(S^{c})$. The kernel of this map consists of powers of the twist
$f$. Also, the image of this map is of finite index subgroup
containing in the pure subgroup $\Gamma'(m)$ of ${\rm Mod}(S^{c})$.
Hence, the center of this image is trivial by
Lemma\,\ref{Lem-trivial}. Hence, $C(C_{\Gamma'}(f^{n}))$ consists of
powers of $f$ and so, it is isomorphic to $\mathbb{Z}$.

Now, let us verify the second condition. If $f$ is a power of a Dehn
twist about a separating simple closed curve $c$, then the curve $c$
separates $S$ in two surfaces with holes. One of the components is
either a nonorientable surface of genus at least three or an
orientable surface of genus at least one. Let $N_1$ denote this
component. Clearly, $N_1$ contains two-sided simple closed curves
$d$ and $e$ so that $d$ and $e$ intersect transversally one point.
Obviously, $t_{d}^{n}$, $t_{e}^{n}$ are in $C_{{\rm M}(S)}(f)$ for
any integer $n$. Hence, if $t_{d}^{n}$, $t_{e}^{n}$ are in
$\Gamma'$, $t_{d}^{n}$, $t_{e}^{n}$ are in $C_{\Gamma'}(f)$. On the
other hand, if $n>0$ then the Dehn twists $t_{d}^{n}$ and
$t_{e}^{n}$ don't commute, and hence $t_{d}^{n}$, $t_{e}^{n}$ are
not in $C(C_{\Gamma'}(f))$. So, we obtain second condition above.

Similar arguments work for nonseparating simple closed curves with
orientable or nonorientable complement. Indeed, if $f$ is a power of
a Dehn twist about a nonseparating simple closed curve $c$ with
nonorientable or orientable complement, then $S^{c}$ is a
nonorientable surface of genus $g-2 \geq 3$ with holes and two
boundary components or $S^{c}$ is an orientable surface of genus
$\frac{g-2}{2} \geq 1$ with holes and two boundary components,
respectively. Since $g-2 \geq 3$  or $\frac{g-2}{2} \geq 1$,  there
are nontrivial simple closed curves $d$  and $e$ on $S^{c}$ meeting
transversally at one point. Finally, the same argument used in the
above paragraph finishes the proof.
\end{proof}

As mentioned in the introduction, algebraic characterizations of Dehn twists
about nonseparating simple closed curves curves in both orientable
(Theorem 2.1. of \cite{I1}) and nonorientable surfaces (provided that the curve
is not characteristic, Theorem 3.1. of \cite{A}; see also \cite{ASzep};) are already done.
On the other hand, the situation for separating curves is more involved.
First we remark that Theorem 3.2 and Theorem 3.3 of \cite{A} are not correct as they are stated.
E. Irmak informed us about the following counter example constructed by L. Paris.

\begin{figure}[hbt]
 \begin{center}
 \includegraphics[width=11cm]{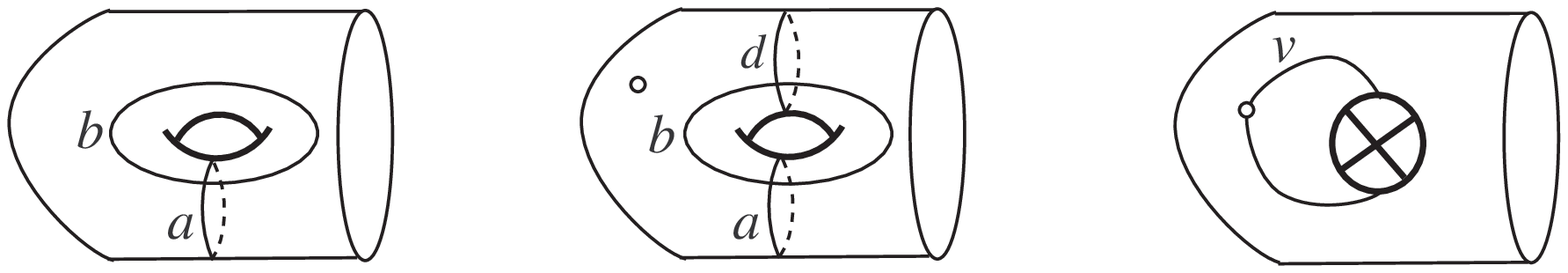}
\caption {} \label{oneholedsurfaces}
\end{center}
\end{figure}

{\bf A Counter Example.} Let \ $c$ \ be a separating curve in a
surface $S$ so that one of the two components of $S^c$ is a torus
with one boundary component $\Sigma_{1,1}$ as in the
Figure~\ref{oneholedsurfaces}. Consider the mapping class
$\tau=(t_at_b)^3$, where $a$ and $b$ are the curves given in the
same figure.  Let $K$ be an abelian subgroup as in the statement of
Theorem 3.2 or Theorem 3.3 of \cite{A} for the Dehn twist $t_c$.
Since $\tau$ \ commutes with both $t_a$ and $t_b$, $\tau$ is in
$C_{{\rm M}(S)}(K)$. However, $t_c=\tau^2$ \ and hence $t_c$ is not
a primitive element in $C_{{\rm M}(S)}(K)$. Note that this example
exists since the center of the mapping class group of $\Sigma_{1,1}$
is not trivial.  Indeed it is isomorphic to the infinite cyclic
group generated by $\tau$. There are two more cases that might cause
similar problem. The first one is $\Sigma_{1,1}^1$, \ whose mapping
class group is isomorphic to that of $\Sigma_{1,1}$. The final
problematic case is the surface $N_{1,1}^1$ \ whose mapping class
group is the infinite cyclic group generated by the class of $v$,
the class of the puncture slide diffeomorphism
(Figure~\ref{oneholedsurfaces}). In this case, $v^2$ \ is the Dehn
twist about the boundary curve $c$.

\begin{Remark}\label{IvnovThm2.2}
The example we have just described above, also provides a counter
example for Theorem 2.2 of \cite{I1}, a result for algebraic
characterization for Dehn twists (mainly, about separating curves)
in orientable surfaces. The methods we will provide here, which work
for both orientable and nonorientable surfaces, not only provide an
algebraic characterization for Dehn twists about separating curves
but also for the topological type of the separating curve the Dehn
twist is about (see Section 5).
\end{Remark}

\section{Preparation for the algebraic characterization of Dehn
twist about separating curves}

We start with the following technical result which we will use
later.

\begin{Proposition}\label{Prop-MaxAbelSub}
Let $S=N_g^k$ be a nonorientable surface of even genus. Then for any
integer $s=0,\cdots, \frac{g-2}{2}$, the group ${\rm M}(S)$ has an
abelian subgroup of rank $\displaystyle\frac{3g-6-2s}{2} + k$, which
is freely generated by Dehn twists about pairwise nonisotopic
nonseparating simple closed curves, and so that no abelian subgroup
containing this subgroup has bigger rank. Moreover, when we cut the
surface along these curves, the resulting surface is a disjoint
union of $g+k-2s-2$ many pair of pants and $2s$ many two holed real
projective planes.
\end{Proposition}

\begin{proof}
Let $r=\displaystyle\frac{3g-6-2s}{2} + k$. By the maximality
condition of the proposition, all the components of the surface cut
along these \ $r$ \ simple closed curves have Euler characteristic
$-1$. In other words, each component is either a pair of pants or a
two holed real projective plane.  Assume that there are $l$ many
pair of pants and $m$ many two holed real projective planes. Hence,
considering Euler characteristics we obtain the equation
$$2-g-k=\chi (S)=\chi (\coprod_l N_0^3 \cup \coprod_m N_1^2) =-l-m \ .$$
On the other hand, counting the number of punctures, we obtain
$$3g-6-2s+3k= 3l+2m \ .$$ These two equations yield $m=2s$ and
$l=g+k-2s-2$ as desired. Finally, the existence of such subgroups is
readily seen by inspection.
\end{proof}

\subsection{Separating chains and pairs of Dehn twists}
A sequence of Dehn twists $t_{a_1},\cdots,t_{a_n},$ is called a
chain if the following geometric intersection $i(a_i,a_{i+1})=1$,
for $i=1,\cdots,n-1$. The integer $n\geq 1$ is called the length of
the chain. Note that if a chain has length more than one then each
$a_i$ must be two-sided nonseparating and has nonorientable
complement. For a tree or chain of Dehn twists we always fix an
orientation for a tubular neighborhood of the union of simple closed
curves (which is always orientable) and consider Dehn twists using
this orientation. It is known that two Dehn twists $t_a, \ t_b$
satisfy the braid relation \ $t_at_bt_a=t_bt_at_b$ \ if and only if
$i(a,b)=1$ on nonorientable surfaces (see \cite{SM1}). Hence, by the
above results any automorphism $\Psi:{\rm M}(S)\rightarrow {\rm
M}(S)$ maps a chain of Dehn twists of length at least two to another
chain of Dehn twists of the same length.  In this note, unless we
state otherwise a chain or a tree in a nonorientable surface $S$
will mean a chain or a tree of Dehn twists about nonseparating
two-sided simple closed curves with nonorientable complement.

Below we will give a generalization of Lemma 3.7 of \cite{A} to
punctured surfaces. We will include the proof since the one in
\cite{A} has a gap, indicated by B. Szepietowski.

\begin{Lemma}\label{Lem-SepChains}
Let $S=N_g^k$ be a nonorientable surface of genus $g\geq 5$ with
$k$. Then the image of a disc separating chain under an automorphism
of ${\rm M}(S)$ is again a chain which separates a disc.
\end{Lemma}

\begin{proof} If the genus of the surface is odd then a chain is
separating if and only if it is maximal in ${\rm M}(S)$.  However,
being maximal is clearly preserved under an automorphism. Now Lemma
3.5 of \cite{A} finishes the proof.

Now assume that the genus $g\geq 6$ is an even integer. Further
assume that $c_1,\cdots, c_{2l+1}$ is a disc separating chain in
${\rm M}(S)$. Hence, when we delete a tubular neighborhood of the
chain from the surface we obtain a disjoint union of a disc and a
nonorientable surface, call $S_0$, of genus $g-2l$ with $k$
punctures and one boundary component. By Euler characteristic
calculation and inspection we see that the group ${\rm M}(S_0)$ has
a abelian subgroup $K$ of rank $\displaystyle\frac{3(g-2l)-6}{2} +
k+2$, contained in each $C_{{\rm M}(S)}(t_{c_i})$, which is freely
generated by Dehn twists about pairwise nonisotopic simple closed
curves.

Now suppose on the contrary that the image $d_1,\cdots, d_{2l+1}$ of
the chain $c_1,\cdots, c_{2l+1}$ under an automorphism is not
separating. So by Lemma 3.5 of \cite{A} the complement of a tubular
neighborhood of the chain $d_1,\cdots, d_{2l+1}$ in $S$ is an
orientable surface of genus $\displaystyle\frac{g-2l}{2}-1$ with $k$
punctures and two boundary components, say $c_1$ and $c_2$. Let us
call this surface $S_1$. (Since the surface $S$ is nonorientable the
curves $c_1$ and $c_2$ are both characteristic in the surface $S$.)
The image of the abelian subgroup $K$ under the same automorphism is
again a maximal subgroup in ${\rm M}(S)$ and it lies in each
$C_{{\rm M}(S)}(t_{d_i})$, $i=1,\cdots, 2l+1$. However, the
orientable subsurface $S_1$ can support an abelian subgroup $K_0$ in
${\rm M}(S)$, which lies in each $C_{{\rm M}(S)}(t_{d_i})$, of rank
at most $\mbox{rank}(K)-1$. Some of the generators of both groups
are Dehn twists about characteristic or separating curves. However,
by Lemma~\ref{Lem-1.1} some powers of these generators are preserved
under automorphisms. This finishes the proof.
\end{proof}

Now we will characterize algebraically a separating pair of Dehn
twists about some two-sided simple closed curves, each of which is
nonseparating with nonorientable complement (so that together they
separate the surface).

\begin{Lemma}\label{Separating Pair Characterization}
Let $S=N_g^k$ be a nonorientable surface of genus $g\geq 5$ and
$a_1$ and $a_2$ be disjoint, nonisotopic, nonseparating two-sided
simple closed curves with nonorientable complements. Then $a_1$ and
$a_2$ together separate the surface if and only if the following
conditions are satisfied:
\begin{enumerate}
\item For any Dehn twist $t_b$ satisfying the braid relation
$t_{a_1}t_bt_{a_1}=t_bt_{a_1}t_b$ we have $t_b\not\in C_{{\rm
M}(S)}(t_{a_2})$;
\item In the even genus case, the twists $t_{a_1}$ and $t_{a_2}$ are
contained in a free generating set, whose elements are all two-sided
nonseparating simple closed curves with nonorientable complement, of
a maximal abelian subgroup $K$ in ${\rm M}(S)$, of rank
$r=\frac{3g-6-2s}{2}+k$, where $s=1$ or $s=2$.
\end{enumerate}
Moreover, if the conditions of lemma are satisfied, then $s=2$ if
both components of the surface $S$ cut along the curves $a_1$ and
$a_2$ are nonorientable of even genus, and $s=1$ otherwise.
\end{Lemma}

\begin{proof}
Suppose first that $a_1$ and $a_2$ separate the surface. Now if a
Dehn twist $t_b$ satisfies the braid relation
$t_{a_1}t_bt_{a_1}=t_bt_{a_1}t_b$ then the curves $b$ and $a_1$
intersects geometrically once. However, since  $a_1$ and $a_2$
separate the surface $b$ and $a_2$ must intersect nontrivially and thus
$t_b$ cannot be contained in $C_{{\rm M}(S)}(t_{a_2})$.
Moreover, it is easy to construct the required free abelian
subgroup $K$, in case the genus of $S$ is even.

For the other direction assume that the conditions of the lemma are
satisfied, but on the contrary suppose that the surface cut along
the curves $a_1$ and $a_2$ is connected.  Hence, the surface $S$ cut
along only $a_2$, say $S_0$, is connected.

First we will treat the case where the genus $g$ is an odd integer.
Hence, the surface $S_0$ is a connected nonorientable surface of
genus at least 3. Moreover, the curve $a_1$ is still nonseparating
in $S_0$. Hence, there is a two-sided curve $b$ in $S_0$, whose
geometric intersection with $a_1$ is one. This is a contradiction to
the first condition of the assumption.

Now let us consider the even genus case.  By
Proposition~\ref{Prop-MaxAbelSub} the surface cut along the $r$ many
curves, about which the Dehn twists generate the subgroup $K$, has
nonorientable components. Thus the surface $S_0$ is a connected
nonorientable surface of genus at least four, and the curve $a_1$ is
neither separating nor characteristic in $S_0$. Hence, as above
there is a two-sided curve $b$ in $S_0$, whose geometric
intersection with $a_1$ is one. This finishes the proof for the even
genus case.

The final part of the lemma is an immediate consequence of
Proposition~\ref{Prop-MaxAbelSub}.
\end{proof}

\subsection{Triangles of Dehn twist}
Let $a_1,a_2,a_3$ be distinct, nonisotopic, nonseparating two-sided
simple closed curves with nonorientable complements. We say that
they form a triangle if each geometric intersection $i(a_i,a_j)=1$,
for all $i\neq j$ (see \cite{ASzep}).  A triangle is called
orientation reversing if there are Dehn twists about these curves,
denoted by $t_{a_1}$, $t_{a_2}$ and $t_{a_3}$, so that
$t_{a_1}t_{a_2}t_{a_1}=t_{a_2}t_{a_1}t_{a_2}$,
$t_{a_2}t_{a_3}t_{a_2}=t_{a_3}t_{a_2}t_{a_3}$ and
$t^{-1}_{a_1}t_{a_3}t^{-1}_{a_1}=t_{a_3}t^{-1}_{a_1}t_{3_2}$. (Note
that on a nonorientable surface we have exactly two Dehn twists
about any two-sided circle $a$, which we may denote by $t_a$ and
$t_a^{-1}$.) Otherwise the triangle is called orientation
preserving. Note that a triangle is either orientation preserving or
orientation reversing but not both (cf. see Theorem 3.1 in
\cite{I1}).  Similarly, we will call the Dehn twists about these
curves orientation reversing or preserving, respectively. It is
clear that this property of triangles of curves or the Dehn twists
about these curves is algebraic and thus it is preserved by
automorphisms of the mapping class groups.

The following topological characterization is proved in \cite{ASzep}.

\begin{Lemma}\label{Punctured Annulus}
The above triangle of Dehn twists is orientation preserving if and
only if the union of these curves has an orientable regular
neighborhood. Moreover, a regular neighborhood of a nonorientable
triangle is $N_{4,1}$, a genus four nonorientable surface with one
boundary component.
\end{Lemma}

\subsection{Dehn twists about characteristics curves and separating curves}
Let $t_a$  be a Dehn twist about a separating simple closed curve
$a$ in a nonorientable surface $S$. Maximal chains contained in the
centralizer $C_{{\rm M}(S)}(t_a)$ correspond to maximal chains of
Dehn twists about two-sided simple closed curves (nonseparating with
nonorientable complement) lying in different components of $S^a$.
The lengths of these maximal chains determine the topological type
the curve $a$ up to a great extend, however fail to characterize its
topological type completely. A more powerful tool is to use maximal
trees contained in the centralizers. By a maximal tree of Dehn
twists in a surface nonorientable $S$ we will mean a connected
maximal tree of Dehn twists about pairwise nonisotopic nonseparating
two-sided simple closed curves with nonorientable complements. If
$S$ is orientable we will only require that the curves in the tree
to be nonseparating.

\begin{Remark}\label{Remark-MaxChainInTree vs In Surface}
As it is seen in the figure below a maximal chain in a maximal tree
need not to be a maximal chain in the surface. Note that the chain
of circles $1,2,\cdots,7$ is maximal in the tree but not in the
surface, which contains the longer chain $1,2,\cdots,7,8,9$.

\begin{figure}[hbt]
 \begin{center}
 \includegraphics[width=6cm]{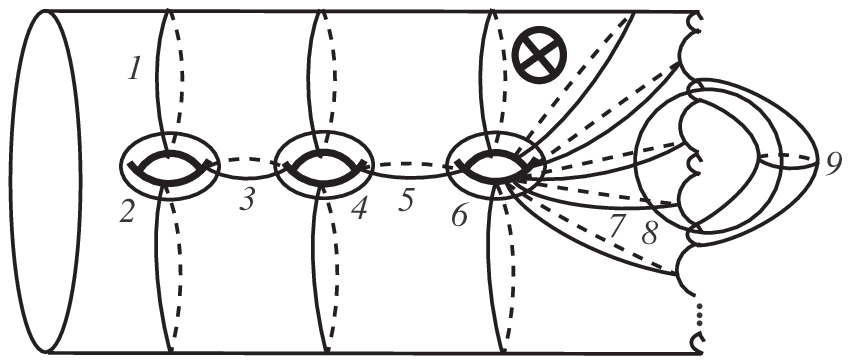}
\caption {} \label{maxtree}
\end{center}
\end{figure}
\end{Remark}

The trees below will be useful for the rest of the paper:
$T_{2g+1,1}^k$,  $T_{2g+2,1}^k$ and $OT_{g,1}^k$.

\begin{figure}[hbt]
 \begin{center}
 \includegraphics[width=8cm]{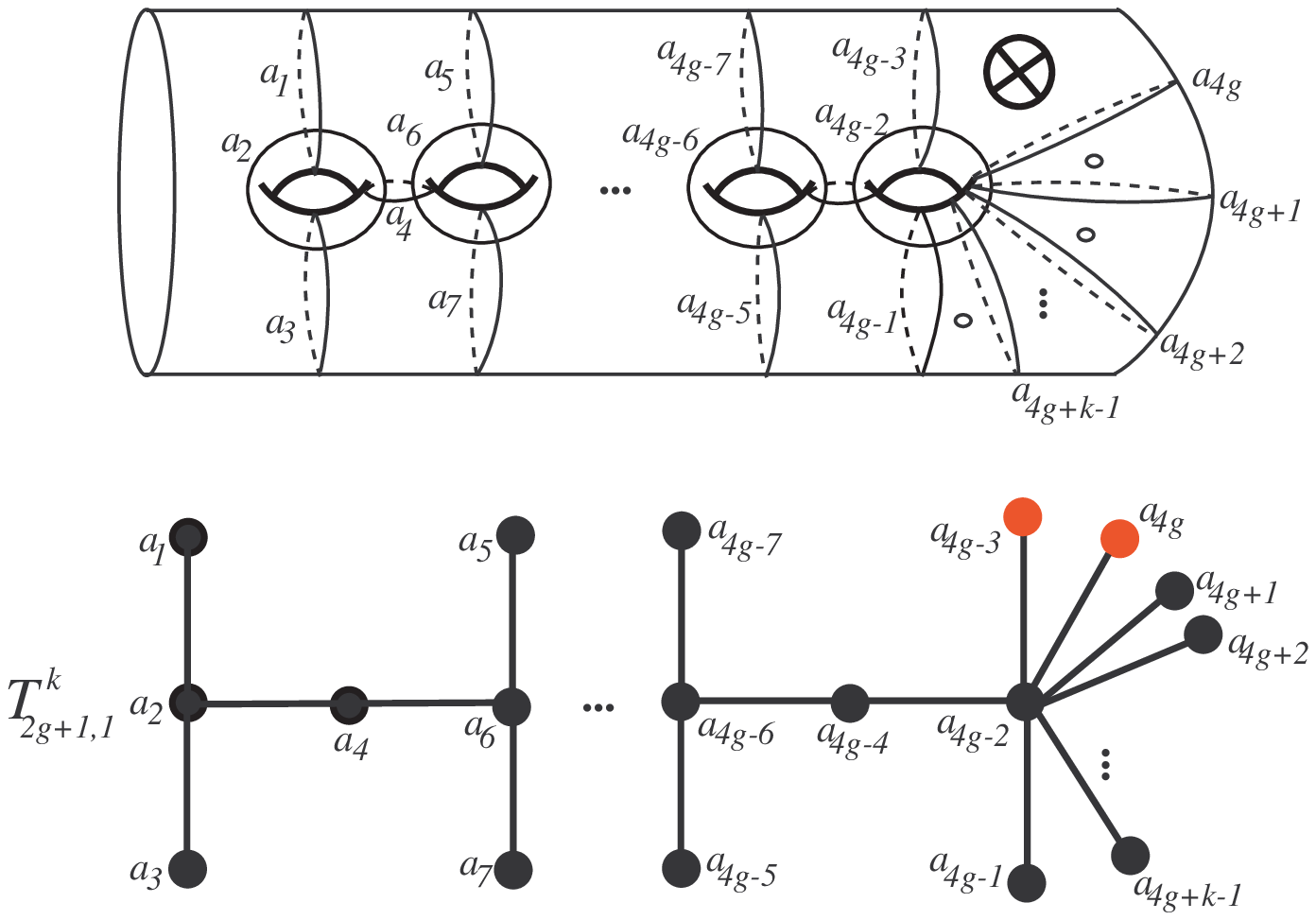}
\caption {} \label{oneholedtree1}
\end{center}
\end{figure}

\begin{figure}[hbt]
 \begin{center}
 \includegraphics[width=8cm]{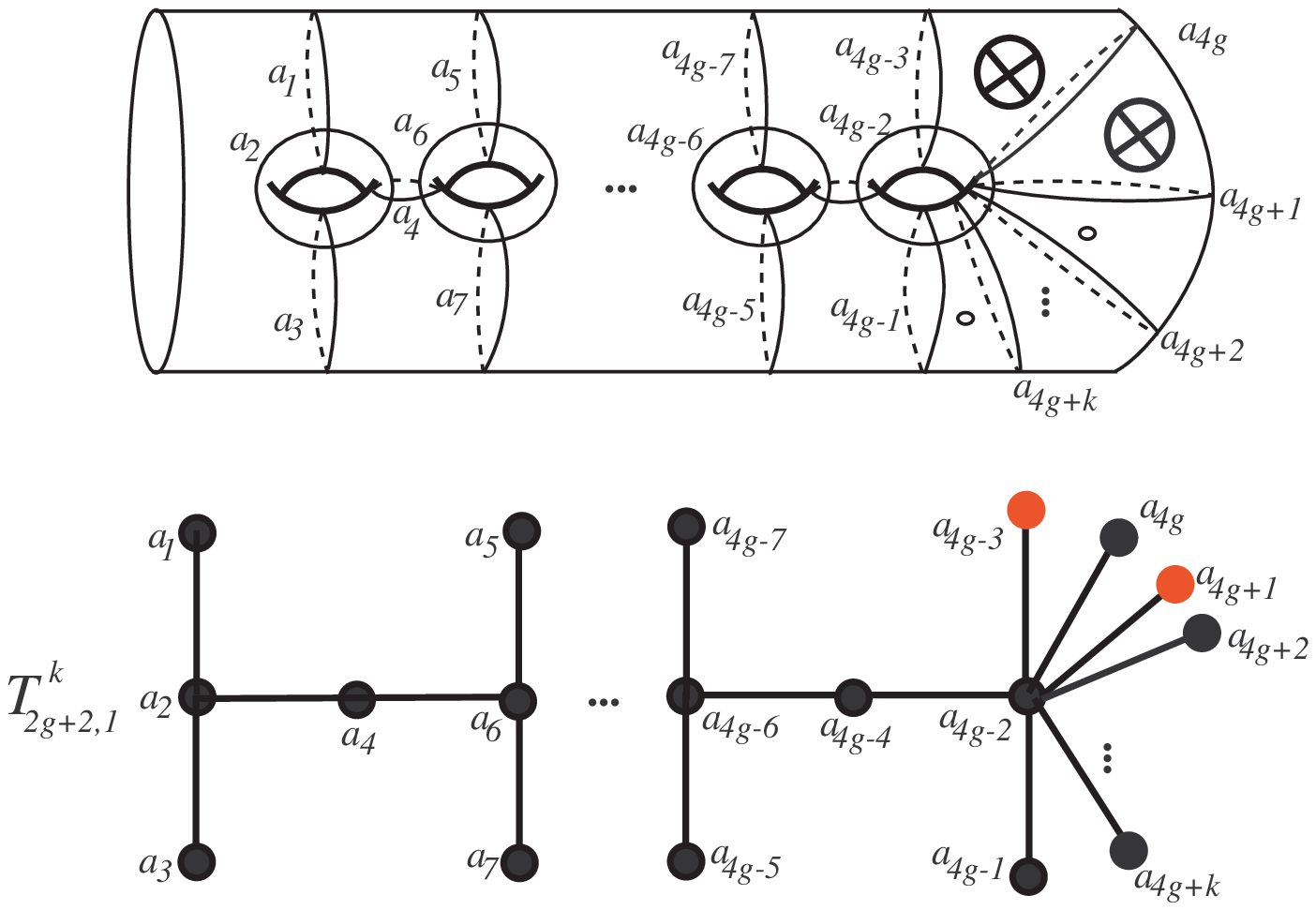}
\caption {} \label{oneholedtree2}
\end{center}
\end{figure}

\begin{figure}[hbt]
 \begin{center}
 \includegraphics[width=8cm]{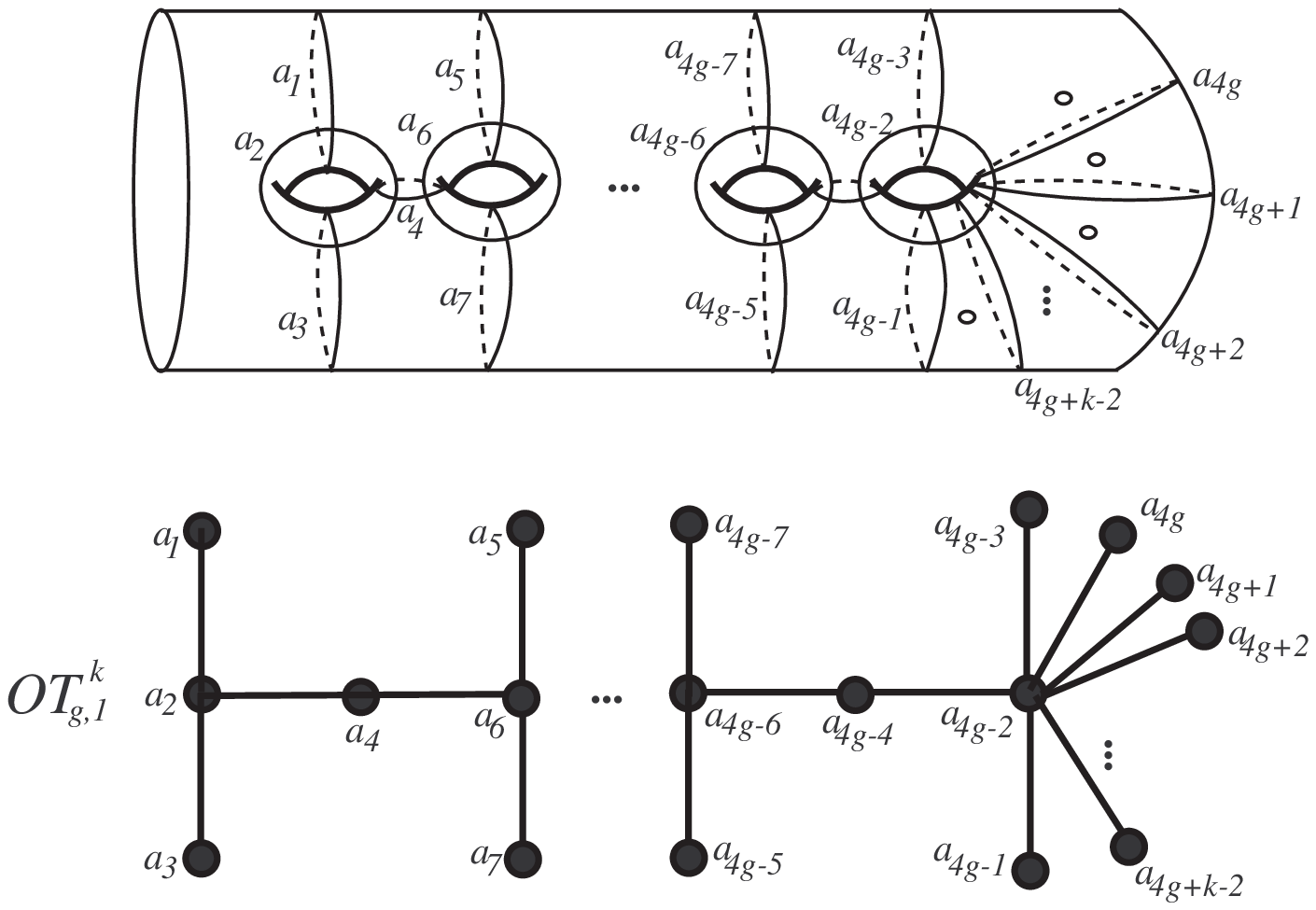}
\caption {} \label{oneholedtree}
\end{center}
\end{figure}

We may endow an embedding of one of these trees into the group ${\rm
M}(S)$, where $S$ is a nonorientable surface of genus at least five,
with a coloring of its vertices.  For example, the tree
$T_{2g+1,1}^k$ in Figure~\ref{oneholedtree1}  has the vertices
$a_{4g-3}$ and $a_{4g}$ colored. This will mean that any orientation
reversing triangle in ${\rm M}(S)$, whose vertices commute with the
colored vertices also commute with all the vertices of the tree.
Equivalently, any orientation reversing triangle which lies in
$$C_{{\rm M}(S)}(t_{a_{4g-3}})\cap C_{{\rm M}(S)}(t_{a_{4g}}) \ ,$$
also lies in $C_{{\rm M}(S)}(t_{a})$, for all vertices $a$ of the
tree.

For nonorientable surfaces of even genera, first we present an
algebraic characterization of a Dehn twist about a two-sided simple
closed curve, whose complement is orientable.

\begin{Lemma}\label{Lem-compOrien}
Let $g\geq 2, \ k\geq 0$ be integers and $c$ be a nontrivial
two-sided simple closed curve in $S=N_{2g+2}^k$. Then $S^c$ is
orientable if and only if the colored tree $T=NT_{2g+2}^k$ (see
Figure\,\ref{tree}) can be embedded in the centralizer $C_{{\rm
M}(S)}(t_c)$ as a maximal tree, where
\begin{enumerate}
\item each maximal chain in the tree is a maximal chain in ${\rm M}(S)$;

\item the tree $T$ has a chain with length larger than or equal to any
chain in the centralizer $C_{{\rm M}(S)}(t_c)$;

\item any two vertices connected to $a_{4g-3}$, both different than
$a_{4g-4}$, form a separating pair.
\end{enumerate}

\end{Lemma}

\begin{proof}
One direction is clear. Now assume that $NT_{2g+2}^{k}$ lies in the
centralizer $C_{{\rm M}(S)}(t_c)$. We claim that the surface $S$ cut
along the curves $a_1$, $a_2$ and $a_3$ has two components one of
which is an orientable surface of genus $g-1$ with $k$ punctures and
one boundary component: To see this, consider a tubular neighborhood
of the tree $NT_{2g+2}^k$ with the vertex $a_0$ deleted. This is an
orientable surface of genus $g$ with $1+k+2(g-1)$ boundary
components.  By the maximality of the tree, $2g-2$ many of these
components must bound discs. To illustrate this, consider, for
example, the maximal chain $a_1,a_2,a_4,a_6,a_5$. The boundary
component corresponding to this chain should bound a disc, because
otherwise the tree would not be maximal (we could attach another
two-sided simple closed curve to $a_2$). Note also that by the
condition (3) of the hypothesis of the lemma the $k$ pairs
$(t_{a_{4g-2}},t_{a_{4g-1}})$, $(t_{a_{4g-1}},t_{a_{4g}})$,
$\cdots$, $(t_{a_{4g+k-3}},t_{a_{4g+k-2}})$ on the right corner of
the tree are all separating. By maximality and the coloring they
must all bound punctured annuli. This finishes the proof of the
claim.

Now, if we attach a tubular neighborhood of the circle $a_0$ to this
subsurface we obtain another subsurface, call $S_0$, of genus $g$
with $k$ punctures and two boundary components. Note that the two
boundary components of $S_0$ should be glued so that the resulting
surface would be nonorientable of genus $2g+2$ with $k$ punctures.
Finally, $c$ being disjoint from each vertex of the tree implies
that $c$ lies in $S-\cup_{i\geq1} a_i$, \ which is a disjoint union
of $2g$ discs, $k$ punctured annuli and a one holed Klein bottle. In
particular, up to isotopy, $c$ is the unique nontrivial two-sided
simple closed curve in this holed Klein bottle.  This finishes the
proof.

\begin{figure}[hbt]
 \begin{center}
 \includegraphics[width=8cm]{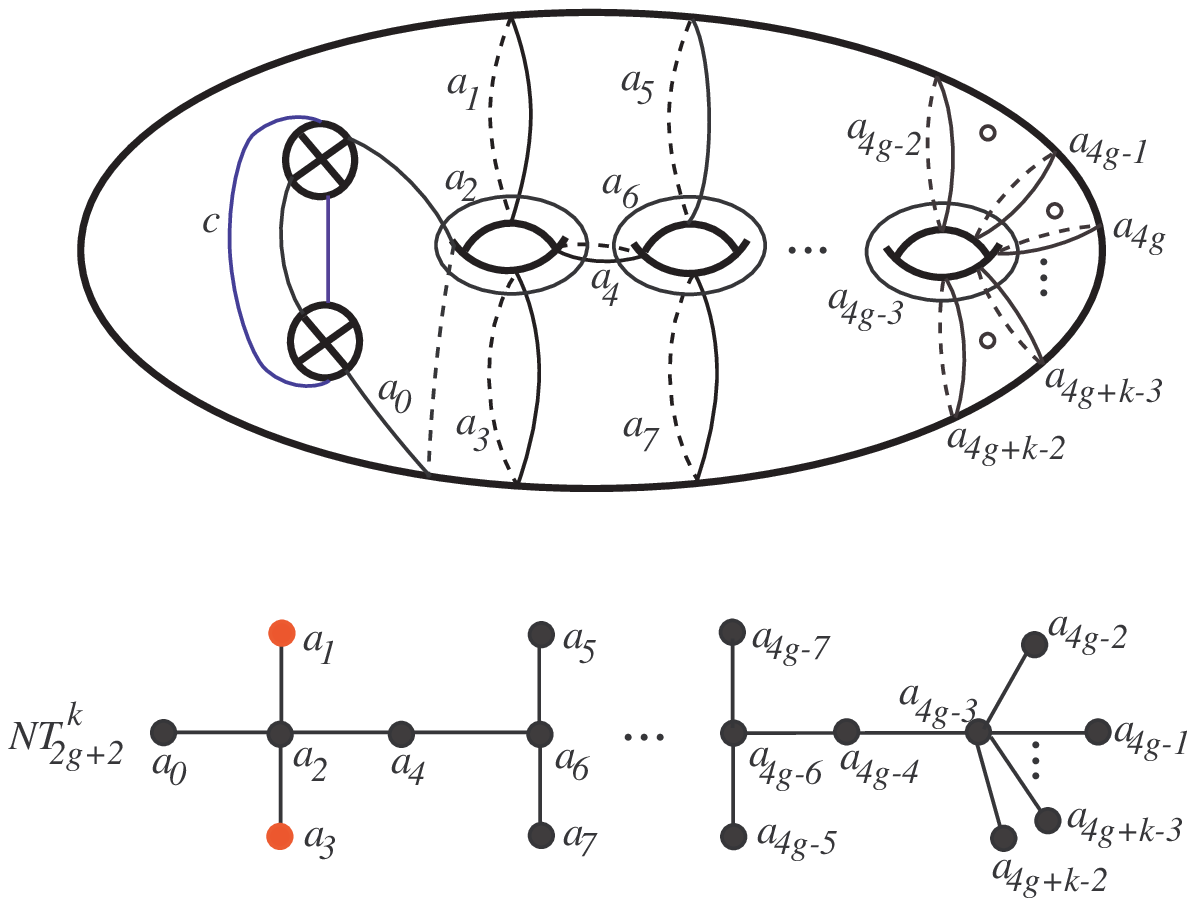}
\caption {} \label{tree}
\end{center}
\end{figure}
\end{proof}

We note that the above lemma and Lemma~\ref{Lem-1.1}
yield the following algebraic characterization
for powers of Dehn twists about nonseparating simple closed curves with
orientable complement on nonorientable surfaces of even genus.

\begin{Corollary}\label{Cor-CompOrient}
Let $g\geq 2, \ k\geq 0$ be integers $f$ a mapping class in ${\rm
M}(S)$ such that $f=t_c^m$ for some integer $m>0$, and a nontrivial
simple closed curve $c$ on $S=N_{2g+2}^k$. Then $c$ is a
characteristic curve (i.e., its complement $S^c$ is orientable) if
and only if the colored tree $T=NT_{2g+2}^k$ (see
Figure\,\ref{tree}) can be embedded in the centralizer $C_{{\rm
M}(S)}(t_c)$ as a maximal tree, where
\begin{enumerate}
\item each maximal chain in the tree is a maximal chain in ${\rm M}(S)$;

\item the tree $T$ has a chain with length larger than or equal to any
chain in the centralizer $C_{{\rm M}(S)}(t_c)$;

\item any two vertices connected to $a_{4g-3}$, both different than
$a_{4g-4}$, form a separating pair.
\end{enumerate}
\end{Corollary}

\bigskip
\noindent

For separating Dehn twists characterization we will make use of the
following lemma.

\begin{Lemma}\label{LemChr-2.2} Let $g$ be a positive integer and
$T$ be one of the colored trees $T_{2g+1,1}^k$, $T_{2g+2,1}^k$ or
$OT_{g,1}^k$, embedded in the group ${\rm M}(S)$, where $S$ is a
nonorientable surface of genus at least five. Suppose that $c$ is a
nontrivial separating simple closed curve in $S$ and the tree $T$
lies in the centralizer $C_{{\rm M}(S)}(t_c)$ as a maximal tree.
Moreover, assume the followings:

\begin{enumerate}
\item Each maximal chain in the tree is a maximal chain in ${\rm M}(S)$;

\item If the tree $T$ is \ $T_{2g+1,1}^k$ \ or \ $T_{2g+2,1}^k$, then it
has a chain with length larger than or equal to any chain in the
centralizer $C_{{\rm M}(S)}(t_c)$.

\item Any two vertices connected to $a_{4g-2}$, except $a_{4g-4}$, form a
separating pair.

\item[(4a)] If $T=T_{2g+1,1}^k$ then $g\geq 2$, and if $T=T_{2g+2,1}^k$ then
$g\geq 1$. Moreover, in both cases, there is an orientation
reversing triangle in $C_{{\rm M}(S)}(t_c)$, which is not contained
in $$\bigcap_{a_i\in V(T)}C_{{\rm M}(S)}(t_{a_i}),$$ where $V(T)$ is
the set of vertices of $T$;

\item[(4b)] If $T=OT_{g,1}^k$ then $g\geq 2$, and any orientation reversing
triangle in ${\rm M}(S)$ also lies in
$$C_{{\rm M}(S)}(t_{a_{4g-3}})\cap C_{{\rm M}(S)}(t_{a_{4g-1}}) \
.$$
\end{enumerate}

\noindent Then $S^c$ has a component homeomorphic to $N_{2g+1,1}^k$,
$N_{2g+2,1}^k$ or $\Sigma_{g,1}^k$, respectively.
\end{Lemma}

\begin{proof}
\noindent {\it Case 1: $T=T_{2g+1,1}^k$}. First we assume that $S$
is of odd genus. For each Dehn twist belonging the tree choose a
two-sided nonseparating simple closed curve $a_i$ with nonorientable
complement.  Let $S_0$ be the closure of a tubular neighborhood of
the tree of the curves $a_i$. Then $S_0$ is an orientable subsurface
of $S$ with Euler characteristic $\chi(S_0)=2-4g-k$ and with $2g+k$
boundary components. So $S_0$ is an orientable surface of genus $g$
with $2g+k$ boundary components. By the definition of maximal tree
each maximal chain in $T$ is a maximal chain in the surface $S$.
Moreover, $S$ has odd genus and thus each maximal chain contained in
$T$ separates the surface (this is the only place we use the
assumption that $S$ is of odd genus). Hence, each boundary component
bounds either a disc, a once punctured disc, an annulus or a
M\"obius band. Again by maximality of the tree, $2g-2$ many of these
components must bound discs. To illustrate this, consider for
example, the maximal chain $a_1,a_2,a_4,a_6,a_5$. The boundary
component corresponding to this chain should bound a disc, because
otherwise the tree would not be maximal (we could attach another
two-sided simple closed curve to $a_2$).

Note that by the condition (3) of the hypothesis of the lemma the
$k+1$ pairs $(t_{a_{4g-3}},t_{a_{4g}})$,
$(t_{a_{4g}},t_{a_{4g+1}})$, $\cdots$,
$(t_{a_{4g+k-1}},t_{a_{4g-1}})$ on the right corner of the tree are
all separating.

The condition (2) of the hypothesis of the lemma implies that at
most two of these pairs may bound a projective plane with two
boundary components. On the other hand, the condition (4a) implies
that at least one of them bounds a projective plane with two
boundary components.  Finally, the coloring of the vertices implies
that boundary components corresponding to the chain
$a_{4g-3},a_{4g}$ is the only pair that bounds a projective plane
with two boundary components. Hence, the other $k$ separating pairs
must bound punctured annuli.

By attaching $2g-2$ discs and $k$ punctured discs and a M\"obius
band to $S_0$ we get a nonorientable surface, say $S_1$, of genus
$2g+1$ with one boundary component. $S_1$ is contained as a
subsurface in one of the two components of $S^c=S_2 \cup S_3$, say
in $S_2$. Since $S_2$ has only one boundary component, the boundary
component corresponding to the chain \ $a_1,a_2,a_3$ \ must be
parallel to the boundary of $S_2$. Hence, we are done in the odd
genus case.

Now let us consider the case where $S$ has even genus. If each
maximal chain in $T$ is separating in $S$ then the above proof works
in this case as well. Now we will show that any maximal chain in $T$
is indeed separating in $S$. To prove this, assume that there is a
maximal chain $t_{c_1},\cdots,t_{c_{2l+1}}$ in $T$, and thus in
${\rm M}(S)$, so that the chain of two-sided simple closed curves
$c_1,\cdots,c_{2l+1}$ is not separating in $S$. A tubular
neighborhood $\nu$ of the union of these curves is an orientable
surface of genus $l$ with two boundary components. Since the chain
is maximal in ${\rm M}(S)$ we see that the surface $S \setminus
int(\nu)$ is an annulus, possibly with punctures, so that when the
orientable surfaces $\nu$ and $S \setminus int(\nu)$ are glued along
the two boundary components, the resulting surface $S$ is
nonorientable (see Figure\,\ref{2holedsurface}) with genus $2l+2$.
Since the separating curve $c$ is disjoint from the tree and thus
from the chain, we see that $c$ lies in the annulus and bounds a (at
least twice) punctured disc.

\begin{figure}[hbt]
 \begin{center}
 \includegraphics[width=8cm]{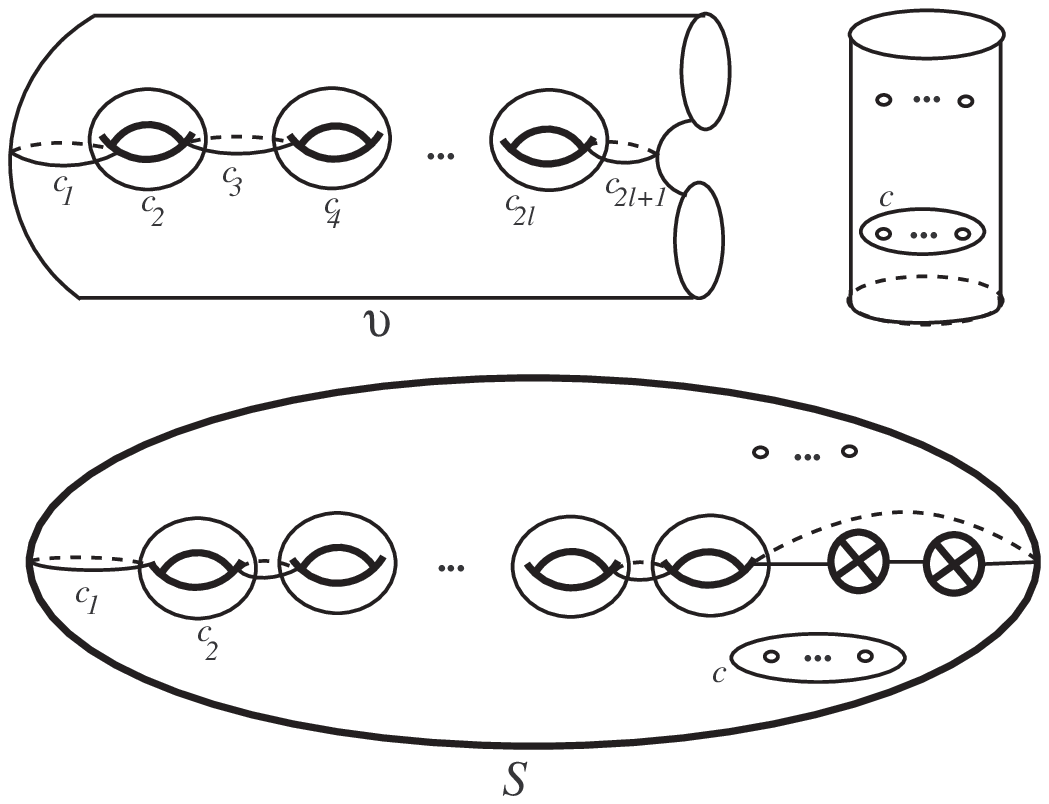}
\caption {} \label{2holedsurface}
\end{center}
\end{figure}

So the topological type of $c$ is determined up to the number of
punctures of $S$ contained in either sides of $c$. On the other
hand, since the surface $S$ has even genus, exactly two of the
separating pairs  \ $(t_{a_{4g-3}},t_{a_{4g}})$,
$(t_{a_{4g}},t_{a_{4g+1}})$, $\cdots$,
$(t_{a_{4g+k-1}},t_{a_{4g-1}})$ bound a projective plane with two
boundary components. On the other hand, the coloring of the vertices
implies that $a_{4g-3},a_{4g}$ is the only pair bound a projective
plane with two boundary components. This gives the desired
contradiction. Hence the proof finishes in the case where
$T=T_{2g+1,1}^k$.
\bigskip

\noindent {\it Case 2: $T=T_{2g+2,1}^k$}. Now first assume that the
surface $S$ has odd genus. We proceed analogous to the previous case
and arrive at the point, where the $k+2$ pairs
$(t_{a_{4g-3}},t_{a_{4g}})$, $(t_{a_{4g}},t_{a_{4g+1}})$, $\cdots$,
$(t_{a_{4g+k}},t_{a_{4g-1}})$ on the right corner of the tree are
all separating. Again the condition (2) of the hypothesis of the
lemma implies that at most two of these pairs bound a projective
plane with two boundary components and the others will bound
punctured annuli. Then condition (4a) and the coloring of the
vertices imply that exactly the two pairs
$(t_{a_{4g-3}},t_{a_{4g}})$ and $(t_{a_{4g}},t_{a_{4g+1}})$ bound
projective plane with two boundary components and all other pairs
bound punctured annuli. This finishes the proof for the odd genus
case.

The even genus case is again analogous to that in Case $1$. We just
need to show that any maximal chain in the tree is separating. We
proceed as in case $T=T_{2g+1,1}^k$. Without loss of generality we
may assume that $c_{2l+1}$ belongs to the set
$$\{a_{4g-3},a_{4g},a_{4g+1},\cdots,a_{4g+k},a_{4g-1}\} \ .$$
The condition (2) of the hypothesis of the lemma implies that the
$k+2$ pairs $(t_{a_{4g-3}},t_{a_{4g}})$,
$(t_{a_{4g}},t_{a_{4g+1}})$, $\cdots$, $(t_{a_{4g+k}},t_{a_{4g-1}})$
on the right corner of the tree are all bounding punctured annuli.
However, this contradicts to the maximality of the tree, since in
this case we may add two more vertices to the tree as in
Figure\,\ref{oneholedtree2}. So we are done in this case too.

\bigskip
\noindent{\it Case 3: $T=OT_{g,1}^k$}. This case is easy now,
because by condition (4b) all the separating pairs on the right side
of the tree will bound punctured annuli (see the paragraph below
Remark~\ref{Remark-MaxChainInTree vs In Surface}). It follows that
the surface $S_1$, in this case, will be an orientable surface of
genus $g$ with one boundary component and with $k$ punctures
(compare with the surface $S_1$ we obtained in the case
$T=T_{2g+1,1}^k$ above).
\end{proof}

\begin{Remark} {\bf 1)} One needs to be careful when using the
above lemma to characterize the topological type of the curve (or
the Dehn twist $t_c$) algebraically: Namely, if one of the
components of $S^c$, the surface $S$ cut along the curve $c$, is
orientable of genus at least two then we can embed $T=OT_{g,1}^k$
into the centralizer $C_{{\rm M}(S)}(t_c)$ satisfying the conditions
of the lemma so that the topological type of the curve (or the Dehn
twist $t_c$) is algebraically characterized. If the orientable
component is of genus one, then the intersection of $C_{{\rm
M}(S)}(t_c)$ with the centralizer of the tree (the intersection of
all the centralizers of the vertices of $T$) contains a pair of Dehn
twists $t_a$ and $t_b$, about nonseparating curves with
nonorientable complements so that $t_a$ and $t_b$ satisfy the braid
relation.

On the other hand, if both components of $S^c$ are nonorientable,
then we have to choose the component of $S^c$ so that condition (2)
of the lemma is satisfied, which is always possible if $g\geq 7$
(see also the next theorem). In other words, we need to choose the
component which contains the longer chain.

{\bf 2)} For odd genus surfaces $S$ the condition (3) of the
hypothesis of the lemma is void. For example consider  the pair
$(t_{a_{4g-3}},t_{a_{4g}})$ in the first tree $T_{2g+1,1}^k$. The
maximal chain  $t_{a_{4g-3}},t_{a_{4g-2}},t_{a_{4g}}$ in $T$ is
maximal in the surface and thus is separating. Hence, a tubular
neighborhood of the chain $t_{a_{4g-3}},t_{a_{4g-2}},t_{a_{4g}}$ is
a torus with two boundary components, each of which is a separating
curve. This implies that the pair $(t_{a_{4g-3}},t_{a_{4g}})$ is
separating.
\end{Remark}

\section{Completing the proof of the algebraic characterization of
Dehn twist about separating curves}

We note that Dehn  twists about two-sided nonseparating simple
closed curves with nonorientable complements is already known
(\cite{ASzep} or \cite{A}). (See Theorem 2.1 in \cite{I1} for
characterization for Dehn twists about nonseparating simple closed
curves on orientable surfaces.) Moreover, by Lemma~\ref{Lem-1.1} and
Corollary~\ref{Cor-CompOrient} it is enough to characterize Dehn
twist about separating curves algebraically, assuming already that
the element have the form $f=t_c^m$, for some integer $m>0$ and a
nontrivial separating simple closed curve $c$ on $S=N_g^k$.  In this
case, one can easily see that there exists a free abelian subgroup
$K$ of ${\rm M}(S)$ generated by $f$ and
$\displaystyle\frac{3g-9}{2}+k$, when $g$ is odd, (respectively,
$\displaystyle\frac{3g-10}{2}+k$, when $g$ is even), twists about
two-sided nonseparating simple closed curves with nonorientable
complements such that $rank \ (K)=\displaystyle\frac{3g-7}{2} +k$,
when $g$ is odd, (respectively, $\displaystyle\frac{3g-8}{2}+k$,
when $g$ is even).

Also note that if genus $g\geq 7$ then one of the two components of $S^c$ is either
a nonorientable surface of genus at least genus four or an
orientable surface of genus at least two. Hence, the above lemma and
the remark following it are applicable.

\begin{Theorem}\label{Chr-2.1} Let $g\geq 5, \ k\geq 0$ be integers $f$ a
mapping class in ${\rm M}(S)$ such that $f=t_c^m$ and $K$ be as above.
If Lemma~\ref{LemChr-2.2} is applicable, which is
always the case if $g\geq 7$, then $f=t_c$ if and only if
\begin{enumerate}
\item $f$ is a primitive element of $C_{{\rm M}(S)}(K)$ if $c$ does not separate
$\Sigma_{1,1}$ \ $\Sigma_{1,1}^1$ \ or $N_{1,1}^1$;
\item $f=(t_at_b)^6$ if $c$ separates $\Sigma_{1,1}$, where $t_a, t_b$ is a chain
contained in the intersection of $C_{{\rm M}(S)}(t_c)$ with the centralizer of the
tree which is used in Lemma~\ref{LemChr-2.2}  (see Figure\,\ref{oneholedsurfaces});
\item $f=(t_at_bt_d)^4$ if $c$ separates $\Sigma_{1,1}^1$, where $t_a,t_b, t_d$
is a chain of Dehn twists, about nonseparating simple closed curves
with nonorientable complements, where all are contained in the
intersection of $C_{{\rm M}(S)}(t_c)$ with the centralizer of the
tree;
\item $f=v^2$ if $c$ separates $N_{1,1}^1$, where $v$ is
the class of the puncture slide diffeomorphism, a generator of the
intersection of $C_{{\rm M}(S)}(t_c)$ with the centralizer of the tree.
\end{enumerate}
Moreover, the topological type of $c$ is determined completely via
Lemma~\ref{LemChr-2.2}.

If $Lemma~\ref{LemChr-2.2}$ is not applicable and $g=6$ then $f=t_c$
if and only if $f$ is primitive in  $C_{{\rm M}(S)}(K)$.
Furthermore, each component of $S^c$ is a nonorientable surface of
genus three and the number punctures in each component is $r_i-2$,
where $r_i$, $i=1,2$ are determined by the maximal trees in
Figure~\ref{daisy} below can be embedded in $C_{{\rm M}(S)}(t_c)$.

Finally, if $Lemma~\ref{LemChr-2.2}$ is not applicable and
$g=5$ then one of the components is again a nonorientable surface of
genus three. The topological type of the other component is also
algebraically characterized.
\end{Theorem}

\begin{proof}
First assume that Lemma~\ref{LemChr-2.2} is applicable. If $c$ does not separate
$\Sigma_{1,1}$ \ $\Sigma_{1,1}^1$ \ or $N_{1,1}^1$ then we are done because the
arguments provided in the proof of Theorem 2.2 of \cite{I1} will
work in this case since the problematic cases are ruled out.

If $c$ separates $\Sigma_{1,1}$ then $f$ is a Dehn twist (i.e.,
$f=t_c$) if and only if $f=(t_at_b)^6$ for some Dehn twists $t_a$
and $t_b$, about nonseparating simple closed curves with
nonorientable complements so that $t_a, t_b$ is a chain (i.e., they
satisfy the braid relation), where both are contained in the
intersection of $C_{{\rm M}(S)}(t_c)$ with the centralizer of the
tree which is used in Lemma~\ref{LemChr-2.2} (see
Figure~\ref{oneholedsurfaces}).

Similarly, if $c$ separates $\Sigma_{1,1}^1$ then $f$ is a Dehn
twist if and only if $f=(t_at_bt_d)^4$ for a chain of Dehn twists
$t_a,t_b, t_d$, about nonseparating simple closed curves with
nonorientable complements, where all are contained in the
intersection of $C_{{\rm M}(S)}(t_c)$ with the centralizer of the
tree (see Figure~\ref{oneholedsurfaces}).

Finally, if $c$ separates $N_{1,1}^1$
then $f$ is a Dehn twist if and only if $f=v^{\pm 2}$, where $v$ is
the class of the puncture slide diffeomorphism, again contained in the intersection of $C_{{\rm
M}(S)}(t_c)$ with the centralizer of the tree (we know that $t_c=v^2$; see Figure~\ref{oneholedsurfaces}).
This provides an algebraic characterization because $v$ is a generator for the mapping
class group of $N_{1,1}^1$, which is isomorphic to the intersection of $C_{{\rm
M}(S)}(t_c)$ with the centralizer of the tree.

If $g=5$ or $6$ and the hypothesis of
Lemma~\ref{LemChr-2.2} is still satisfied then again there is
nothing to do. Hence, we are left with the cases $g=5$ or $6$ and Lemma~\ref{LemChr-2.2} is
not applicable. Hence, one of the components of $S^c$ is nonorientable of genus three.

In case $g=6$, first note that, both components of $S^c$ are
nonorientable of genus three and thus $f=t_c$ if and only if $f$ is
primitive in  $C_{{\rm M}(S)}(K)$.  Moreover, via Euler
characteristic considerations, the maximal trees in
Figure~\ref{daisy} contained in $C_{{\rm M}(S)}(t_c)$ determines the
number of punctures in each component.  Indeed, it is $r-2$ if $r$
is as in the theorem.

\begin{figure}[hbt]
 \begin{center}
 \includegraphics[width=4cm]{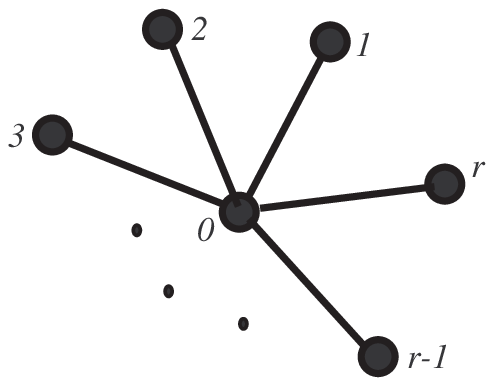}
\caption {} \label{daisy}
\end{center}
\end{figure}

If $g=5$ we proceed as follows: We know that one component is a punctured
nonorientable surface of genus three.  Thus we need a way of checking algebraically
whether the second component is a puncured torus or a Klein bottle.
A punctured nonorientable surface of genus three supports a maximal chain of length three.
A (punctured) torus supports a chain of length at least two and a punctured Klein bottle supports
no chain of length two.  Hence, if $C_{{\rm M}(S)}(t_c)$ contains two maximal chains
of lengths at least two, where any Dehn twist from one chain commutes with any of the
Dehn twists of the other chain, then the second component of $S^c$ is a punctured torus.
Otherwise the other component is a punctured Klein bottle. The number of punctures in each
component can be determined as in the genus six case. Note that the Dehn twists in
Figure~\ref{daisy} contained in one component should commute with all the Dehn twists with
the chain in the other component.  This way we can determine the number of punctures in each component.

If the second component is a punctured Klein bottle or a punctured
torus with more than one punctures then $f=t_c$ if and only if $f$
is primitive in  $C_{{\rm M}(S)}(K)$. If the second component is
torus with at most one puncture then we proceed as in the parts (2)
and (3) above.
\end{proof}

\section{Orientable Surfaces}\label{OrientSurf}
In this section, analogously to the above section, we will give an
algebraic characterization for Dehn twists about separating curves on
an orientable surface of genus at least three.  We will omit the
proofs since they are essentially the same. In fact, the versions for
orientable surfaces will be even easier since a large portion of our
efforts has been spent to distinguish algebraically an annulus with
one puncture from the projective plane with two boundary components.
First let us give a version of Lemma~\ref{LemChr-2.2} for orientable
surfaces.

\begin{Lemma}\label{LemChr-5.1} Let $g$ be a positive integer and
$T$ be the colored tree $OT_{g,1}^k$, embedded in the group ${\rm
M}(S)$, where $S$ is an orientable surface of genus at least two.
Suppose that $c$ is a nontrivial separating simple closed curve in
$S$ and the tree $T$ lies in the centralizer $C_{{\rm M}(S)}(t_c)$
as a maximal tree. Moreover, assume the followings:
\begin{enumerate}
\item Each maximal chain in the tree is a maximal chain in ${\rm M}(S)$;

\item Any two vertices connected to $a_{4g-2}$, except $a_{4g-4}$, form a
separating pair.
\end{enumerate}
\noindent Then $S^c$ has a component homeomorphic to
$\Sigma_{g,1}^k$.
\end{Lemma}

Now we can state the characterization result for separating Dehn twists:

\begin{Theorem}\label{Chr-5.1} Let $g\geq 3, \ k\geq 0$ be integers $f$ a
mapping class in ${\rm M}(S)$ such that $f=t_c^m$ for some
integer $m>0$ and a nontrivial separating simple closed curve
$c$ on the orientable surface $S=\Sigma_g^k$. If Lemma~\ref{LemChr-5.1} is applicable, which is
always the case if $g\geq 4$, then $f=t_c$ if and only if
\begin{enumerate}
\item $f$ is a primitive element of $C_{{\rm M}(S)}(K)$ if $c$ does not separate
$\Sigma_{1,1}$ or $\Sigma_{1,1}^1$;
\item $f=(t_at_b)^6$ if $c$ separates $\Sigma_{1,1}$, where $t_a, t_b$ is a chain
contained in the intersection of $C_{{\rm M}(S)}(t_c)$ with the centralizer of the
tree which is used in Lemma~\ref{LemChr-5.1}  (see Figure\,\ref{oneholedsurfaces});
\item $f=(t_at_bt_d)^4$ if $c$ separates $\Sigma_{1,1}^1$, where $t_a,t_b, t_d$
is a chain of Dehn twists, about nonseparating simple closed curves,
where all are contained in the intersection of $C_{{\rm M}(S)}(t_c)$
with the centralizer of the tree (see
Figure\,\ref{oneholedsurfaces}).
\end{enumerate}
Moreover, the topological type of $c$ is determined completely via
Lemma~\ref{LemChr-5.1}.

If $Lemma~\ref{LemChr-5.1}$ is not applicable then the curve
separates a punctured torus. The number punctures in each component
can be determined the same way as in Theorem~\ref{Chr-2.1}. If the
number of punctures is at least two then $f=t_c$ if and only if $f$
is primitive.  Otherwise, parts (2) and (3) of this theorem applies.
\end{Theorem}

\section{Some Applications}
Let $S=N_{g}^{k}$. The subgroup $\mathcal{T}$ of the mapping class
group generated by all Dehn twists about two-sided simple closed
curves, is called the twist subgroup. It is known that this subgroup
is of index $2^{k+1} k!$ in ${\rm Mod}(S)$, provided that $g \geq 3$
(\cite{K}, \cite{SM2}). Now, as a consequence of characterization
results of Dehn twists we state the following corollary:

\begin{Corollary}\label{Cor-AlgChr}
For $g\geq 5$, let $\Phi:{\rm M}(S)\rightarrow{\rm M}(S)$ be an
automorphism and $\mathcal{T}\leq {\rm M}(S)$ be the twist subgroup.
If $t_c\in \mathcal{T}$ is a Dehn twist then so is $\Phi(t_c)$.
Moreover, the Dehn twists $t_c$ and $\Phi(t_c)$ are topologically
equivalent. In other words, there is a homeomorphism $f:S
\rightarrow S$ \ such that $\Phi(t_c)=t_{f(c)}$. In particular, the
twist subgroup is a characteristic subgroup of ${\rm M}(S)$.
\end{Corollary}

The subgroup ${\rm PMod}^{+}(S)$ has index $2^k$ in ${\rm PMod}(S)$
and contains the twist subgroup $\mathcal{T}$ as a subgroup of index
two (\cite{SM2}) provided that $g\geq 3$.

\begin{Lemma}\label{Lem-Chr-PMPlus}
The subgroup ${\rm PMod}^{+}(S)$ is characteristic in ${\rm Mod}(S)$
and ${\rm PMod}(S)$, provided that $g\geq 5$.
\end{Lemma}

\begin{proof}
We know that the twist subgroup is characteristic in all the three
groups in the statement of the lemma. On the other hand, ${\rm
PMod}^{+}(S)$ is generated by the twist subgroup and the set of all
$Y$-homeomorphisms, each of which is supported inside Klein bottles
with one boundary component, so that both components of the surface
cut along this boundary curve are nonorientable. It is easy to see
that this set is also characteristic (cf. Theorem 3.9 and Theorem
3.10 of \cite{A}). Indeed, one can see this directly as follows:
Note that, in the Klein bottle with the boundary circle $e$, a
$Y$-homeomorphism represents a mapping class, say $\tau$, which is
characterized as a mapping class that is not a Dehn twist but
$\tau^2=t_e$. It follows that ${\rm PMod}^{+}(S)$ is characteristic
in ${\rm Mod}(S)$ and ${\rm PMod}(S)$.
\end{proof}

\section*{acknowledgement}

I would like to thank B. Szepietowski for his numerous valuable
comments on the earlier version of this paper. I would like to thank
also  Y. Ozan for reviewing and several suggestions.

\bigskip
\providecommand{\bysame}{\leavevmode\hboxto3em{\hrulefill}\thinspace}

\end{document}